# ESTIMATION OF A FUNCTION UNDER SHAPE RESTRICTIONS. APPLICATIONS TO RELIABILITY

### By L. Reboul

*Université de Poitiers and Université de Marne-La-Vallée*

This paper deals with a nonparametric shape respecting estimation method for U-shaped or unimodal functions. A general upper bound for the nonasymptotic $\mathbb{L}_1$-risk of the estimator is given. The method is applied to the shape respecting estimation of several classical functions, among them typical intensity functions encountered in the reliability field. In each case, we derive from our upper bound the spatially adaptive property of our estimator with respect to the $\mathbb{L}_1$-metric: it approximately behaves as the best variable binwidth histogram of the function under estimation.

**1. Introduction.** In this paper we study a data-driven nonparametric estimation method for shape restricted functions. As an application of this study, we first have in mind classical frameworks such as estimation of unimodal densities or regression functions. We also place stress on building and studying shape respecting estimators of typical intensity functions, namely the hazard rate of an absolutely continuous distribution and the failure rate of a nonhomogeneous Poisson process, which are key concerns in systems reliability studies: for a nonreparable system, which is replaced by a new one after it fails, the failure behavior is modeled by the distribution of its single lifetime, frequently specified via its hazard rate. For reparable systems, repaired but not replaced after each failure, the failure behavior in time can be modeled by a counting process. When repair times can be disregarded and the system has a large number of units, this counting process can be approximated by a nonhomogeneous Poisson process. Such a process is totally characterized by its cumulative intensity function or, when it exists, by its failure rate.

Nonparametric estimation procedures have often been investigated first for density estimation and regression and then generalized to other frameworks. The more widely used are smoothing or projection methods with fixed parameters (see [24] and [9] for densities and [25] for regression functions). Several estimators of this type have been proposed and studied for the hazard rate, under censored and uncensored schemes (see [28, 30, 33]). In other respects, Curioni [11] studies histograms and kernel estimators of the failure rate of a nonhomogeneous Poisson process, based on the observation of replications of the process. Even if the difficult problem of the choice of the smoothing parameter can be overcome by automatic methods such as cross-validation (see [26] for density estimation, [19] and [10] in reliability), the main handicap of those methods lies in their rigidity: they tend to assume that the unknown function has homogeneous variation everywhere. In other words, they are not sensitive enough to the local concentration of the data. Such a drawback clearly appears in the problem of estimating the hazard rate by fixed bandwidth kernel estimators: the local variance of those estimators tends to increase towards infinity as the number of systems at risk decreases. Obviously, these methods are totally misleading for estimating the failure rate of a nonhomogeneous Poisson process in realistic situations, where one generally observes a small number of replications of the failure time process on a finite time period. Indeed, the system's condition at time $t$ depends on its whole history before $t$ so that the situation is truly nonasymptotic. One therefore needs estimation methods flexible enough to balance the lack of information collected by fitting the data as well as possible, making a locally sensitive choice of the parameter. For that purpose, variable bandwidth kernel estimators and variable binwidth histograms have been studied, first for densities by Stone [29], and then for the hazard rate by Müller and Wang [21, 22]. Bartoszyński, Brown, McBride and Thompson [3] propose a variable bandwidth kernel estimator of the failure rate, based on the observation of replications of the process. The choice of the local bandwidth is generally done by the minimization of an asymptotic mean square error estimator.

Another way of building adaptive tools is to look for the nonparametric maximum likelihood estimate over a restricted class of functions, under which the likelihood is to be maximized. Contrary to kernel estimators, the construction of these estimators does not require either a smoothing parameter or any smoothness assumption on the unknown function and only relies on very natural shape restrictions. Brunk [8] proposes the isotonic estimator for monotone regression functions and Durot [14] studies its good asymptotic properties related to the $\mathbb{L}_1$-metric. Such estimators for decreasing hazard rates have been put forward by Barlow, Bartholomew, Bremner and Brunk [1] in complete life data models. Similarly, Bartoszyński, Brown, McBride and Thompson [3] and Barlow, Proschan and Sheuer [2] propose the nonparametric maximum likelihood estimator for decreasing failure rates. The shape restriction in the last two cases is very natural since it corresponds to the observation of a system during its debugging period. For a decreasing density, the nonparametric maximum likelihood estimate is known as the Grenander estimator [15]. It has a very simple graphical meaning since it is the slope of the least concave majorant of the empirical distribution function based on a sample generated by the density under estimation. It takes the form of a variable binwidth histogram, generating a partition which is approximately the best one in the $\mathbb{L}_1$-metric sense. This property is checked by Birgé [5, 6] from a nonasymptotic minimax risk point of view and Groeneboom [16] and Groeneboom, Hooghiemstra and Lopuhaä [17] study its good asymptotic $\mathbb{L}_1$-properties. The construction of the Grenander estimator and its properties can straightforwardly be extended to the case of a unimodal density with known mode. Nevertheless, a more realistic assumption is that the mode is unknown. Actually, the nonparametric maximum likelihood estimator does not exist any more on such a wide class. One can solve the problem (see [31, 32]) by finding a prior estimate of the mode, but the resulting properties rely on the choice of this estimate. In several studies, Birgé [4, 7] proposes a totally data-driven estimation method for unimodal densities with unknown mode: his estimator relies neither on the arbitrary choice of extra parameters nor on any smoothness assumption on the unknown density. It still approximately behaves as the best histogram in terms of the nonasymptotic minimax $\mathbb{L}_1$-risk, over restricted sets of unimodal densities.

Our purpose in this paper is to extend Birgé's method to a more general functional estimation framework. More precisely, we have in mind to define and study estimators for positive integrable functions $g$, assumed to be unimodal or U-shaped (decreasing then increasing). The unimodal assumption is often realistic for regression or density functions, while U-shaped hazard rate or failure rate functions correspond to the failure behavior of a system which is observed during its entire lifetime: after a debugging period where the number of failures tends to decrease, the latter is stable during the exploitation period, and then turns out to deteriorate from aging. Starting from a step function estimator $\widehat{G}$ of $G = \int_a^{\cdot} g(t)\, dt$ on $I = [a, b]$, we define the shape respecting estimator $\widetilde{g}$ of $g$ as the image of $\widehat{G}$ through some deterministic mapping. The definition of this mapping relies on a convenient adaptation of the "Pool Adjacent Violators Algorithm" (see [1]) which is involved in the definition of Grenander estimator.

This paper is organized as follows: in Section 2, we define and study this mapping in a deterministic framework. The former study is applied in Section 3 to a statistical framework: we build a general upper bound for the $\mathbb{L}_1$-risk of the shape respecting estimator of U-shaped or unimodal functions and investigate conditions under which $\widetilde{g}$ behaves as a "clever" histogram, generating on its own a partition which is optimal from an $\mathbb{L}_1$-risk point of view. Section 4 is devoted to the application of our results to particular functions: we first study the shape respecting estimator of unimodal regression and density functions. We next build and study shape respecting estimators for a U-shaped hazard rate and a U-shaped failure rate in realistic underinformed designs. The proofs of our results are given in the three last sections.

**2. Deterministic framework.** Our general aim in this paper is to estimate a shape restricted function $g$ defined on a given compact real interval $I = [a, b]$. In Section 3 we will show that we can define a shape respecting estimator of $g$ as the image of a step function estimator of $G = \int_a^{\cdot} g(t) \, dt$, through a deterministic mapping. This section is thus devoted to the construction and the study of such a mapping. More precisely, we are interested in mappings from the cone $\mathcal{H}(I)$ of nondecreasing, right-hand continuous with left-hand limits (cadlag) step functions on $I$ into particular sets of shape restricted, integrable functions on $I$. We focus here on the set of U-shaped functions and the set of unimodal functions on $I$, defined as follows:

DEFINITION 1. A function $g$ defined on an interval $I = [a, b]$ is a U-shaped (resp. unimodal) function if there exists some number $m$ in $I$ such that $g$ is nonincreasing (resp. nondecreasing) on $[a, m]$ and nondecreasing (resp. nonincreasing) on $[m, b]$.

For the sake of simplicity, we shall restrict ourselves to the U-shaped case. The unimodal case is briefly described in Remark 3. Moreover, we approach $G$ by a nondecreasing function, which implies that $g$ is positive, but this assumption may probably be dropped (see Remark 5).

2.1. *Construction of the mapping.* Let $m \in I$ be an arbitrary point. To be clearer, we first define a mapping $U_{\mathcal{G}}^m$ from $\mathcal{H}(I)$ into $\mathcal{U}_{\mathcal{G}}^m(I)$, where $\mathcal{U}_{\mathcal{G}}^m(I)$ is the set of U-shaped functions on $I$ whose minimum is achieved at $m$ (when $m$ coincides with one of the endpoints of $I$, we get the important subsets of nonincreasing and nondecreasing functions). This mapping generalizes the isotonic mapping classically used in various contexts of functional estimation under monotonicity restrictions: when $g$ is a decreasing function on $I = [a, b]$, it is defined as the slope of the least concave majorant of an approximation $F$ of $G$ and can be computed using the "Pool Adjacent Violators Algorithm" (PAVA), described in [1]. Formally, the mapping $U_{\mathcal{G}}^m$ is defined as follows:

DEFINITION 2. Let $I = [a, b]$ be a compact real interval and let $m$ be an arbitrary point in $I$. Let $F \in \mathcal{H}(I)$. We define $U_{\mathcal{G}}^m(F)$ as the right-hand continuous slope of $\widetilde{F}_S^m$, where $\widetilde{F}_S^m$ is defined on $[a,m]$ as the least concave majorant of the restriction of $F$ to $[a,m]$ and is defined on $[m,b]$ as the greatest convex minorant of the restriction of $F$ to $[m,b]$. The function $\widetilde{F}_S^m$ is called the U-shaped regularization at $m$ of $F$ on $I$.

Let us notice that $\widetilde{F}_S^m$ is a continuous, piecewise affine function on $I$. We now turn to the definition of our main mapping $U_S$ from $\mathcal{H}(I)$ into $\mathcal{U}_S(I)$, where $\mathcal{U}_S(I) = \bigcup_{m \in I} \mathcal{U}_S^m(I)$ is the set of all U-shaped functions on $I$. For this purpose, we use an idea introduced by Birgé [4] in the context of the estimation of a unimodal density. It consists in minimizing on $I$ the function $d_S$ defined by

$$d_S(m) = \sup_{t \in I} |F(t) - \widetilde{F}_S^m(t)|.$$

It is easy to see that $d_S$ is a cadlag step function on $I$, whose discontinuity points belong to the set of discontinuity points of $F$. This property gives a sense to the following:

DEFINITION 3. Let $F \in \mathcal{H}(I)$ and let $m(F)$ denote the midpoint of the interval where the function $d_S$ defined above achieves its minimum. We define the mapping $U_S$ from $\mathcal{H}(I)$ into $\mathcal{U}_S(I)$ by $U_S(F) = U_S^{m(F)}(F)$.

It is worth noticing that $U_S(F)$ is easily computable in practice. The determination of $m(F)$ does not cause any trouble since we need to compare only a finite number of regularizations. Moreover, for any $m \in I$, the regularization $\widetilde{F}_S^m$ can be computed via the PAVA: the algorithm is applied on $[a,m]$ to compute the least concave majorant of the restriction of $F$ to $[a,m]$, then on $[m,b]$ to compute the greatest convex minorant of the restriction of $F$ to $[m,b]$.

2.2. *An $\mathbb{L}_1$-approximation upper bound.* We now investigate the $\mathbb{L}_1$-properties of our mapping. To fix ideas, let us take the point of view of approximation theory: let $I = [a,b]$ be a compact real interval and let $F \in \mathcal{H}(I)$ be an approximation of some $G = \int_a^{\cdot} g(t)\,dt$, where $g$ is a positive U-shaped function on $I$. Then $U_S(F)$ is a shape respecting approximation of $g$. We seek to link the $\mathbb{L}_1$-approximation quality of $g$ by $U_S(F)$ to the properties of the underlying error $F - G$. In order to square with our stochastic framework, we will have a more general approach: starting from some $F$ of $\mathcal{H}(I)$, we want to control the $\mathbb{L}_1$-distance between $U_S(F)$ and any U-shaped function $g$.

From now on, we will adopt the following:

NOTATION 1. (i) Let $I = [a, b]$ denote a compact real interval, let $f$ and $g$ belong to $\mathbb{L}_1(I)$ and let $\mathcal{H}$ be a subset of $\mathbb{L}_1(I)$. We set

$$\|f - g\| = \int_I |f(t) - g(t)| \, dt \quad \text{and} \quad d(f, \mathcal{H}) = \inf_{h \in \mathcal{H}} \|f - h\|.$$

(ii) Let $\pi$ be a finite partition of $I$ in intervals such that $\pi = ([t_{k-1}^\pi, t_k^\pi))_{k \in \mathcal{K}^\pi}$. Here, $a = t_0^\pi < \cdots < t_{D^\pi}^\pi = b$, $\mathcal{K}^\pi = \{1, \ldots, D^\pi\}$ and $D^\pi$ is arbitrary. We denote by $\Pi(I)$ the set of all such partitions.

(iii) Let $\pi \in \Pi(I)$. We denote by $\mathcal{H}_\pi$ the set of cadlag step functions based on $\pi$ (we mean cadlag step functions that are constant on every $[t_{k-1}^\pi, t_k^\pi)$).

With the above notation, we can state our main theorem (the proof is postponed to Section 5):

THEOREM 1. *Let $I = [a, b]$ be a compact real interval. Let $g$ be a U-shaped function on $I$ and let $G = \int_a^\cdot g(t) \, dt$. Let $F \in \mathcal{H}(I)$ and let $U_S(F)$ be defined by Definition 3. Setting $Z = F - G$ and $f = U_S(F)$, there exists some absolute constant $C \geq 1$ such that*

$$(1) \quad \|f - g\| \leq \inf_{\pi \in \Pi(I)} \left\{ 4 d(g, \mathcal{H}_\pi) + C \sum_{k \in \mathcal{K}^\pi} \sup_{t \in [t_{k-1}^\pi, t_k^\pi]} |Z(t) - Z(t_{k-1}^\pi)| \right\}.$$

*($C = 49$ works.)*

### 2.3. *Comments.*

REMARK 1. Let us denote by $\mathcal{R}_Z(\pi)$ the term in brackets in (1). The quantity $\mathcal{R}_Z(\pi)$ is the sum of a perturbation term $d(g, \mathcal{H}_\pi)$, relying on the smoothness of $g$, and of a regularization term measuring the approximation error of $G$ both by $F$ and by a U-shaped regularization of $F$. Theorem 1 thus stresses that the step function $U_S(F)$ realizes over $\Pi(I)$ the best compromise between those two terms. Let us notice that the perturbation term varies smoothly with $\pi$, while the regularization term achieves its infimum over the subset of partitions whose endpoints belong to the set of discontinuity points of $F$. Thus, the best trade-off will be reached for such a partition.

REMARK 2. From an approximation theory point of view, such a tool, which is able to make a sensitive choice of the image set $\mathcal{H}_\pi$ among all partitions of $I$, is much more powerful than classical projection operators mapping onto linear sets generated by a uniform partition. Indeed, let $\sigma$ and $\pi$, respectively, denote a uniform and a nonuniform partition of $I$ of cardinality $D$. Let $p_\sigma g$ and $p_\pi g$ be the orthogonal projections of some function $g$ on $\mathcal{H}_\sigma$ and $\mathcal{H}_\pi$, respectively. It is proved (see, e.g., [13]) that to achieve an approximation error $\|g - p_\sigma g\| = \mathcal{O}(D^{-\alpha})$, one needs to impose on $g$ smoothness conditions of the type $w(g) = \mathcal{O}(D^{-\alpha})$ (where $w$ is the continuity modulus of $g$), while for approximation by $p_\pi g$, the conditions are less demanding: smoothness conditions on $g$ are needed only in the larger space $\mathbb{L}_\gamma$, $\gamma = (\alpha + 1)^{-1}$. Typically, for $g \in B^\alpha_{p,p}$ with $p < 1$, $\alpha = 1/p - 1$, the approximation error is of order $D^{-\alpha}$ for $p_\pi g$, while for $p_\sigma g$ one has to impose the condition $g \in B^\alpha_{1,1}$ to achieve this rate.

REMARK 3. Let us briefly describe the unimodal case. Let $F \in \mathcal{H}_\pi$ and let $m \in I = [a, b]$. Symmetrically, we can define the unimodal regularization $\widetilde{F}^m_N$ at $m$ of $F$ on $I$ as follows: on $[a, m]$ it is the greatest convex minorant of the restriction of $F$ to $[a, m]$; on $[m, b]$ it is the least concave majorant of the restriction of $F$ to $[m, b]$. We then define $U^m_N(F)$ as its right-hand continuous slope. As $\widetilde{F}^m_S$, the function $\widetilde{F}^m_N$ is piecewise affine. But while $\widetilde{F}^m_S$ is always continuous, $\widetilde{F}^m_N$ happens to be discontinuous at the point $m$ whenever $F$ is discontinuous at this point. Next, in order to define the shape restricting mapping $U_N$, we minimize on $I$ the function $d_N$ defined by

$$d_N(m) = \sup_{t \in I} |F(t) - \widetilde{F}^m_N(t)|.$$

As shown by Birgé [7], $d_N$ is a continuous function on $I$ and its minimum is achieved at a unique point $m(F)$, which is a continuity point of $F$. We thus define the mapping $U_N$ from $\mathcal{H}(I)$ into the set of unimodal functions on $I$ by $U_N(F) = U^{m(F)}_N(F)$. The practical construction of $U_N(F)$ is done by applying the symmetric procedure with the PAVA, after having found $m(F)$. It may be computed numerically by using, for instance, a dichotomous algorithm, after having found the interval where the minimum of $d_N$ is achieved. Setting $f = U_N(F)$, (1) holds for every unimodal function $g$ on $I$. The proof can be performed by exchanging the roles of the intervals $[a, m]$ and $[m, b]$.

REMARK 4. The compactness assumption on $I$ is made for sake of simplicity. Actually it is sufficient to restrict oneself to intervals $I$ on which $G$ is bounded.

REMARK 5. In this paper the underlying function $F$ is assumed to be nondecreasing. This means that our approximation method is applied in practice to positive functions $g$. On the other hand (see, e.g., the application to the estimation of a regression function), one may wish to estimate a function $g$ that is not positive. In such cases we will make the following conjecture: *at the expense of subsidiary technical complications, the monotonicity restriction on $F$ can be dropped.* Checking this conjecture leads namely to show that the minimum of the function $d_S$ is still well defined.

**3. Statistical framework.** The previous results can be applied in a statistical context to build shape respecting estimators for U-shaped or unimodal functions $g$ defined on a given interval $I = [a, b]$. Let $X$ be a random variable whose law depends on $g$, where $g$ is assumed to belong to $\mathcal{U}_S(I)$ [resp. $\mathcal{U}_N(I)$]. Let $\widehat{G} \in \mathcal{H}(I)$ be an estimator of $G = \int_a^{\cdot} g(t)\, dt$ based on the observation $X$ (e.g., $X$ is a sample generated by an unknown unimodal density $g$ and $\widehat{G}$ is the empirical distribution function of the sample). We can apply the mappings previously defined on $\widehat{G}$ to build a shape respecting estimator $\widetilde{g}$ of $g$.

The aim of this section is to study the nonasymptotic properties of this estimator. We first state in Theorem 2 a stochastic version of Theorem 1 that gives control of the $\mathbb{L}_1$-risk of $\widetilde{g}$. We next investigate some conditions of optimality of this control.

In order to be clearer, we still restrict ourselves to the U-shaped case, although the results still hold in the unimodal case. Moreover, we use Notation 1. As a straightforward consequence of Theorem 1, we get:

THEOREM 2. *Let $g$ be a U-shaped function on $I = [a, b]$ and let $X$ be a random variable whose law depends on $g$. Let $\widehat{G} \in \mathcal{H}(I)$ be an estimator of $G = \int_a^{\cdot} g(t)\, dt$ based on $X$ and let $U_S(\widehat{G})$ be defined by Definition 3. Setting $Z = \widehat{G} - G$ and $\widetilde{g} = U_S(\widehat{G})$, there exists some absolute constant $C \geq 1$ such that*

$$(2) \qquad \mathbb{E}\|\widetilde{g} - g\| \leq \inf_{\pi \in \Pi(I)} \mathcal{R}_Z(\pi),$$

*where*

$$\mathcal{R}_Z(\pi) = 4d(g, \mathcal{H}_\pi) + C \sum_{k \in \mathcal{K}^\pi} \mathbb{E}\left( \sup_{t \in [t_{k-1}^\pi, t_k^\pi]} |Z(t) - Z(t_{k-1}^\pi)| \right).$$

*($C = 49$ works.)*

Theorem 2 emphasizes the adaptive behavior of our tool: without advanced knowledge of $g$, it makes a sensitive choice of the partition of $\Pi(I)$ that minimizes the quantity $\mathcal{R}_Z(\pi)$; this quantity takes the form of a risk (it is the sum of a bias term and of the expectation of a random error term). As $\widetilde{g}$ is by construction a histogram based on a random partition, such a sensitive behavior leads us to wonder about its quality, compared to classical histogram estimators of $g$.

Histogram estimators are common tools for estimating a density function $g$ defined on some interval $I$: given an arbitrary partition $\pi \in \Pi(I)$, the histogram estimator of $g$ based on $\pi$ is the empirical estimator of the orthogonal projection of $g$ on $\mathcal{H}_\pi$. Similarly, given a general function $g$ and a step function estimator $\widehat{G}$ of $G$, we consider as an estimator of $g$ the histogram $\widehat{g}^{\pi}$ defined by

$$(3) \qquad \widehat{g}^{\pi}(t) = \sum_{k \in \mathcal{K}^{\pi}} \frac{\widehat{G}(t_k^{\pi}) - \widehat{G}(t_{k-1}^{\pi})}{t_k^{\pi} - t_{k-1}^{\pi}} \mathbb{1}_{[t_{k-1}^{\pi}, t_k^{\pi})}(t) \qquad \text{for all } t \in I.$$

We now investigate conditions on $\widehat{G}$ under which $\widetilde{g}$ turns out to do at least as well as any variable binwidth histogram $\widehat{g}^{\pi}$ of $g$ built from $\widehat{G}$. For that task we can show the following result (see Section 6 for the proof).

THEOREM 3. *Let $\widehat{g}^{\pi}$ be defined by* (3). *Assume that the conditions of Theorem* 2 *hold. Assume moreover that there exists some positive constant $A$ such that for all $\pi \in \Pi(I)$ and all $k \in \mathcal{K}^{\pi}$,*

$$(4) \qquad \mathbb{E}\left(\sup_{t \in [t_{k-1}^{\pi}, t_k^{\pi}]} |Z(t) - Z(t_{k-1}^{\pi})|\right) \leq A\mathbb{E}(|Z(t_k^{\pi}) - Z(t_{k-1}^{\pi})|).$$

*Then we get for all $\pi \in \Pi(I)$,*

$$(5) \qquad (CA + 8)^{-1} \mathcal{R}_Z(\pi) \leq \mathbb{E}\|\widehat{g}^{\pi} - g\| \leq \mathcal{R}_Z(\pi).$$

*Moreover,*

$$(6) \qquad \mathbb{E}\|\widetilde{g} - g\| \leq (CA + 8) \inf_{\pi \in \Pi(I)} \mathbb{E}\|\widehat{g}^{\pi} - g\|.$$

The quality of estimation appears to rely on the expected variations of the process $Z = \widehat{G} - G$ on $I$: for sufficiently small ones, $\widetilde{g}$ generates on its own a partition which is optimal from an $\mathbb{L}_1$-risk point of view. In particular, it will do better than any variable binwidth histogram. Moreover, it is totally data-driven, which is a serious practical advantage since it allows one to solve the problem of how to check for the best partition. Further applications will show that the condition (4) often holds.

REMARK 6. Theorem 2 gives control of the $\mathbb{L}_1$-risk of the estimator of a U-shaped function $g$ with unknown minimum point. Now, suppose that $g$ is U-shaped and that the location of its minimum $m$ is known. An obvious estimator $\widetilde{g}^m$ of $g$ can be obtained via the basic mapping $U_S^m$ of Definition 2 as $\widetilde{g}^m = U_S^m(\widehat{G})$ (when $g$ is a decreasing density, it is merely a Grenander estimate). An upper bound for the $\mathbb{L}_1$-risk of this estimator is derived in Section 5.1 as a by-product of the proof of Theorem 1. Indeed, it follows from Lemma 1 that there exists some absolute constant $C'$ such that

$$(7) \qquad \begin{aligned} &\mathbb{E}\|\widetilde{g}^m - g\| \\ &\leq \inf_{\pi \in \Pi(I)} \left\{ 4d(g, \mathcal{H}_{\pi}) + C' \sum_{k \in \mathcal{K}^{\pi}} \mathbb{E}\left(\sup_{t \in [t_{k-1}^{\pi}, t_k^{\pi}]} |Z(t) - Z(t_{k-1}^{\pi})|\right) \right\}. \end{aligned}$$

Hence, even if the point $m$ is known, one does not lose much by assuming that it is unknown and (2) is still a good control for the $\mathbb{L}_1$-risk of $\widetilde{g}^m$, from a qualitative point of view.

**4. Applications.** We now present examples of application of our study. We first consider the shape respecting estimation of two classical functions: that of a unimodal density, which allows one to recover Birgé's results, and that of a unimodal regression function. For this task, we use the mapping $U_N$ defined in Remark 3. Our main focus is on the estimation of classical intensity functions used in reliability theory, in realistic designs: for nonreparable systems, we study the shape respecting estimator of a U-shaped hazard rate in right-censoring life data models. In the reparable system field, we study the shape respecting estimator of the U-shaped failure rate of a nonhomogeneous Poisson process, based on the observation of a single process on a finite time period. Such a process is widespread in reparable systems studies, since it models the failure behavior of a system having a large number of units and whose repair times can be disregarded. It is totally characterized by its failure rate. In these two last applications, the U-shaped assumption is very natural and corresponds to the situation where a system is observed during its entire lifetime. Moreover, the adaptivity property is particularly important here since realistic reliability designs are often underinformed.

For each application, $I = [a, b]$ is an interval and $\widehat{G} \in \mathcal{H}(I)$ is an estimator of $G = \int_a^{\cdot} g(t)\,dt$, where $g$ is the function under estimation. Moreover, we adopt Notation 1 and the further:

NOTATION 2. We denote by $\Pi_D(I)$ the subset of $\Pi(I)$ of partitions in $D$ intervals. Next, for a given $\pi \in \Pi(I)$, we call $\widehat{g}^\pi$ the histogram estimator of $g$, defined by (3).

Short proofs of the results presented in the sequel are postponed to Section 7.

4.1. *Estimation of a unimodal density.* Let $(X_1, \ldots, X_n)$ be a sample generated by an absolutely continuous distribution $G$ with density $g$. Assume that the restriction of $g$ to a given real interval $I$ is unimodal ($I$ can be the real line here, see Remark 4). The shape respecting estimator of $g$ on $I$ can be defined by $\widetilde{g} = U_N(\widehat{G})$, where

$$\widehat{G}(t) = \frac{1}{n}\sum_{i=1}^{n} \mathbb{1}_{X_i \leq t} \qquad \text{for all } t \in I.$$

As an application of Theorems 2 and 3, we get:

PROPOSITION 1. *There exists an absolute positive constant $E$ such that*

$$(8) \qquad \mathbb{E}\|\widetilde{g} - g\| \leq \inf_{D \leq n}\left\{\inf_{\pi \in \Pi_D(I)} 4d(g, \mathcal{H}_\pi) + E\sqrt{\frac{D}{n}}\right\}.$$

*Moreover, there exists an absolute positive constant $K$ such that*

$$(9) \qquad \mathbb{E}\|\widetilde{g} - g\| \leq K \inf_{\pi \in \Pi(I)} \mathbb{E}\|\widehat{g}^\pi - g\|.$$

Inequality (8) is similar from a qualitative point of view to Birgé's Theorem 1 [7] and gives an idea of how the shape respecting estimator operates: it first chooses among the $D$-dimensional linear subsets of step functions the one which is closest to the unknown function $g$; it then checks the dimension $D$ which realizes the best trade-off between the bias and the error terms of the estimation. Moreover, (9) shows that the selected partition is optimal from a nonasymptotic $\mathbb{L}_1$-risk point of view. This point had been investigated by Birgé from a minimax point of view, since he showed that his estimator nearly achieves the minimax risk over the class of unimodal densities with bounded support. We get here a result for every unimodal density function.

4.2. *Estimation of a unimodal regression function.* Another classical problem is that of the estimation of a unimodal regression function. Let us consider here the model

$$Y_i = g(x_i) + \varepsilon_i, \qquad i = 1, \ldots, n,$$

where $Y_i$ is the observation at time $x_i = i/n$, the $\varepsilon_i$'s are i.i.d. Gaussian centered errors with variance $\sigma^2 > 0$ and $g$ is a unimodal function on $I = [0, 1]$. We can define the shape respecting estimator of $g$ by $\widetilde{g} = U_N(\widehat{G})$, where

$$\widehat{G}(t) = \frac{1}{n}\sum_{i=1}^n Y_i \mathbb{1}_{x_i \leq t} \qquad \text{for all } t \in [0, 1].$$

Let us notice that $\widehat{G}$ is not nondecreasing in general even if $g$ is positive. We thus need to assume that the conjecture of Remark 5 holds in order to define and study $\widetilde{g}$.

PROPOSITION 2. *Assume that the conjecture of Remark 5 holds. Then there exists some positive constant $E$ that depends only on $\sigma$ and $M = \sup_{t \in [0,1]} g(t)$ such that (8) holds. Moreover, (9) holds for some absolute positive constant $K$.*

REMARK 7. The assumption $x_i = i/n$, $i = 1, \ldots, n$, can be dropped as soon as the $x_i$'s are approximately uniformly spread in $[0, 1]$. Moreover, no particular assumption is required on the distribution family of the errors to get (8), and the common distribution has only to be bounded (not necessarily Gaussian) to get (9).

4.3. *Estimation of a U-shaped hazard rate in right-censoring life data models.* Suppose one observes $n$ copies of a nonreparable system in a usual right-censoring scheme: let $(T_1, \ldots, T_n)$ be the potential lifetimes of the $n$ copies, generated by an absolutely continuous distribution $F$ with density $f$. Let $(U_1, \ldots, U_n)$ be the sample of times of censure (this means that $U_j$ is the time beyond which the $j$th copy can no longer be observed). Assume moreover that the $U_j$'s are independent of the $T_j$'s. The random variable $T_j$ will be observed whenever $T_j \leq U_j$. The set of observable data is thus given by

$$\{X_j = T_j \wedge U_j, \delta_j = \mathbb{1}_{X_j = T_j}, j = 1, \ldots, n\}.$$

We want to estimate the hazard rate $g = f/(1 - F)$ of the system's lifetime $T$ on some compact real interval $I = [0, c]$, under the assumption that $g$ is U-shaped on $I$. A shape respecting estimator $\widetilde{g} = U_S(\widehat{G})$ of $g$ on $I$ can be derived from the Nelson–Aalen estimator [23] of the log-survival function $G$. It is defined on $I$ by

$$\widehat{G} = \int_0^\cdot \frac{d\widehat{F}(s)}{1 - \widehat{F}^-(s)}.$$

Here $\widehat{F}$ is the Kaplan–Meier product-limit estimator of $F$ [18] and $\widehat{F}^-$ is its left-hand continuous version. We get:

PROPOSITION 3.  *Assume that $L(c) < 1$, where $L$ is the common distribution function of the $X_j$'s. Then there exists an absolute constant $C$ such that* (8) *holds with $E = C(1 - F(c))^{-1/2}$.*

The sensitive behavior of our tool is a serious advantage here. Indeed, one of the typical drawbacks of hazard rate estimation methods is that the error term of the estimation depends on a penalization factor relying on the value of $L$ near the right-hand side point $c$ of the estimation interval, in such a way that for classical methods such as fixed bandwidth kernels or fixed binwidth histograms, the local error term near this point tends toward infinity. Using $\widetilde{g}$ minimizes the problem. Here the penalization factor is $(1 - L(c))^{-1/2}$, but the partition chosen by the estimate will probably lead to a larger local binwidth when estimating near $c$ than when estimating in the interior. Moreover, the automatic choice of the histogram binwidths is a serious practical advantage and its good behavior holds even for moderate sample sizes. For all these reasons, our estimate should be preferred—at least when the hazard rate is not known to be smooth—to data-adaptive bandwidth kernel estimators, for which the choice of the optimal local bandwidths relies on heavy procedures of asymptotic minimization of a mean square error estimate.

REMARK 8. The definition and properties of the shape respecting estimator of $g$ in complete life data models can be straightforwardly derived from the former results, setting $L = F$.

### 4.4. Estimation of the U-shaped failure rate of a nonhomogeneous Poisson process.

Let $(N(t))_{t \geq 0}$ be a nonhomogeneous Poisson process with mean function $\mathbb{E}(N(t)) = G(t)$. Suppose that $(N(t))_{t \geq 0}$ describes the number of failures in time of a reparable system. The failure rate of the process is, when it exists, the derivative $g$ of $G$.

In the sequel, we propose to estimate $g$ on a given finite time period $I = [0, T]$, where it is known to be U-shaped. For this task, we observe the failure times $(T_1, \ldots, T_{N(T)})$ falling in $I$ of a single copy of the system. Let us define on $I$ $\tilde{g} = U_S(\widehat{G})$ and

$$\widehat{G}(t) = \sum_{k=1}^{N(T)} \mathbb{1}_{T_k \leq t} = N(t).$$

In order to describe the $\mathbb{L}_1$-behavior of $\tilde{g}$ on $I$, we need to use a normalized $\mathbb{L}_1$-distance (which allows one to recover the $\mathbb{L}_1$-distance between two constant rates on $I$). We thus define for all $f, g \in \mathbb{L}_1(I)$

$$\|f - g\| = \frac{1}{T} \int_0^T |f(t) - g(t)| \, dt.$$

We can show that:

PROPOSITION 4. *There exists an absolute positive constant $B$ such that*

(10) $$\mathbb{E}\|\tilde{g} - g\| \leq \inf_{D \leq T} \left\{ \inf_{\pi \in \Pi_D(I)} 4d(g, \mathcal{H}_\pi) + B\sqrt{\frac{G(T)}{T}}\sqrt{\frac{D}{T}} \right\}.$$

*Moreover, there exists an absolute positive constant $K$ such that* (9) *holds.*

The main difference between failure rate estimation and the former studied estimation functions is that the underlying observations $(T_1, \ldots, T_{N(T)})$ are not i.i.d., except in the case where $g$ is constant. Therefore, in the realistic context of the observation of a single system (or a small number of copies of the system), nonadaptive methods are totally misleading. Now let us investigate in some particular cases the quality of $\tilde{g}$ with regard to classical estimators: in the case where $g \equiv \lambda$, (10) allows one to recover the parametric rate. Indeed, the parametric maximum likelihood estimate of $g$ is given in this case by $\widehat{g} = N(T)/T$ and satisfies $\mathbb{E}\|\widehat{g} - \lambda\| \leq \sqrt{\lambda/T}$, while our bound gives the same, up to a multiplicative constant. (More generally, we find rates of order $T^{-1/2}$ whenever $g$ is a step function.) Moreover, when the penalization factor $\sqrt{G(T)/T}$ involved in (10) is bounded, which includes nonincreasing and low nondecreasing functions $g$, the order of magnitude of the asymptotic risk will be at least $T^{-1/3}$ as for densities. More generally, from a nonasymptotic point of view, we hope to obtain a good estimator as soon as the penalization factor is not too large. In some cases, this factor can be very important, in such a way that the quality of the estimate can be quite bad, for instance for high increasing rates. (Fortunately, we are not interested in such situations since they practically correspond to a substantial deterioration of the system, which is retrieved as soon as possible from exploitation.) Nevertheless, even in unfavorable situations, the locally sensitive property of $\widetilde{g}$ will allow one to check break points in the slope of $g$. This fact is interesting for its own sake for trend studies.

REMARK 9. An adaptive estimator of $g$ has been proposed by Barlow, Proschan and Scheuer [2] for decreasing failure rates. It is defined as the nonparametric maximum likelihood estimator over this shape restricted class. Although it has been used successfully in practice (it behaves reasonably well compared with classical parametric models used in industrial reliability studies), its properties have not been investigated so far. The shape respecting estimator we present and study in the sequel generalizes Barlow et al.'s estimator to U-shaped failure rates. Since a decreasing $g$ can be seen as a degenerate U-shaped function, with minimum point at the end of the estimation interval, it is worth noticing that our study also applies to their estimate (see Remark 6).

**5. Proof of Theorem 1.** Let $m$ be a point of $I$ where $g$ achieves its minimum. Let $\widetilde{F}^m$ denote the U-shaped regularization of $F$ at $m$ (see Definition 2) and set $f^m = U_S^m(F)$. We shall prove the following two lemmas in Sections 5.1 and 5.2, respectively (we use Notation 1).

LEMMA 1. *There exists an absolute constant* $C' \geq 1$ *such that for all* $\pi \in \Pi(I)$,

$$\|f^m - g\| \leq 4d(g, \mathcal{H}_\pi) + C' \sum_{k \in \mathcal{K}^\pi} \sup_{t \in [t_{k-1}^\pi, t_k^\pi]} |Z(t) - Z(t_{k-1})|. \tag{11}$$

LEMMA 2. *We have*

$$\|f^m - f\| \leq 4 \sup_{t \in I} |Z(t)|.$$

For all $\pi \in \Pi(I)$ and all $k$, we thus get

$$\sup_{t \in [t_{k-1}^\pi, t_k^\pi]} |Z(t)| \leq \sum_{i=1}^k \sup_{t \in [t_{i-1}^\pi, t_i^\pi]} |Z(t) - Z(t_{i-1}^\pi)|,$$

so that

$$\sup_{t \in I} |Z(t)| \leq \sum_{k \in \mathcal{K}^\pi} \sup_{t \in [t_{k-1}^\pi, t_k^\pi]} |Z(t) - Z(t_{k-1}^\pi)|.$$

Moreover,

$$\|f - g\| \leq \|f^m - g\| + \|f - f^m\|$$

and Theorem 1 follows from Lemma 1 and Lemma 2 with $C = C' + 4$.

5.1. *Proof of Lemma* 1. For the sake of simplicity, we omit the subscript $\pi$ in Notation 1 and we adopt the following extra:

NOTATION 3.   (i) For every subinterval $J$ of $I$ and all functions $f$ and $g$ in $\mathbb{L}_1(I)$, we denote by $l(J)$ the length of $J$ and we set

$$\bar{g}(J) = \frac{1}{l(J)} \int_J g(t)\, dt,$$

$$bg(J) = \int_J |g(t) - \bar{g}(J)|\, dt \quad \text{and} \quad \|f - g\|(J) = \int_J |f(t) - g(t)|\, dt.$$

(ii) $\forall k \in \mathcal{K}$, we set $I_k = [t_{k-1}, t_k]$.

Note first that we may assume without loss of generality that there exists some $j \in \mathcal{K}$ such that $m = t_j$. Indeed, assume that Lemma 1 holds over the restricted class of partitions of $I$ including $m$, for some absolute constant $C''$. Then let $\pi^m$ be a partition of $I$ with $D+1$ endpoints $t_k^m$ such that $t_j^m = m$, and let $\pi$ be the partition of $I$ with $D$ endpoints $t_k$ such that $t_k = t_k^m$ for all $k < j$ and $t_k = t_{k+1}^m$ for all $k \geq j$. For the term in $j$ in (11), we get

$$\sup_{t \in I_j^m} |Z(t) - Z(t_{j-1}^m)| + \sup_{t \in I_{j+1}^m} |Z(t) - Z(t_j^m)|$$

$$\leq 2 \sup_{t \in I_j} |Z(t) - Z(t_{j-1})| + |Z(t_j^m) - Z(t_{j-1})|$$

$$\leq 3 \sup_{t \in I_j} |Z(t) - Z(t_{j-1})|.$$

Moreover,

$$d(g, \mathcal{H}_{\pi^m}) \leq d(g, \mathcal{H}_\pi).$$

Thus, we obtain (11) over the set of all partitions of $I$ by adding constants ($C' = 3C''$).

The main argument to show (11) is Proposition 1 of [5], which is the formal translation of the PAVA. Let us recall this result:

LEMMA 3 (Birgé).  *Suppose that we are given a nondecreasing integrable function $h$ and a nonincreasing integrable function $g$ on some finite interval $J$. Then, using Notation 3,*

$$\|\bar{h}(J) - g\|(J) \leq \|h - g\|(J)$$

(1) Let us first see what is happening on $[a, m]$ in the case where $a < m$. Recall that $t_j = m$. For $k \leq j$, define $\widetilde{F}_k$ as the least concave majorant of the restriction of $F$ on $I_k$. Next, let us define $H_0$ on $[a, m]$ by $H_0(m) = \widetilde{F}_k(m)$ and

$$H_0(t) = \sum_{k \leq j} \widetilde{F}_k(t) \mathbb{1}_{[t_{k-1}, t_k[}(t),$$

and let $h_0$ be the right-hand continuous slope of $H_0$. Then $H_0$ is a piecewise affine, continuous function such that on $[a, m]$, $F \leq H_0 \leq \widetilde{F}^m$. Moreover, $h_0$ is well defined on $[a, m)$ and right-hand continuous. Its discontinuity points belong to $X$, where $X = \{x_0, x_1, \ldots, x_n\}$ is the ordered set obtained as the union of the set $\{t_0, t_1, \ldots, t_j\}$ and that of the discontinuity points of $F$. Now, for $l \geq 1$, let us define $H_l$ and $h_l$ by iterating the following rule (we apply the PAVA to $H_{l-1}$):

(a) If $h_{l-1}$ is nonincreasing on $[a, m)$, then we define $h_l = h_{l-1}$ and $H_l = H_{l-1}$.

(b) If $h_{l-1}$ is not nonincreasing on $[a, m)$, then there exists some $0 \leq i < n$ such that $h_{l-1}(x_{i+1}) > h_{l-1}(x_i)$. Let us define

$$i_- = \min_{0 \leq k < n} \{k : h_{l-1}(x_{k+1}) > h_{l-1}(x_k)\}$$

and

$$i_+ = \min_{i_- < k < n} \{k : h_{l-1}(x_{k+1}) \leq h_{l-1}(x_k)\},$$

where by convention, $\inf \varnothing = m$. Then define $H_l$ such that $H_l = H_{l-1}$ on $[a, x_{i_-}] \cup [x_{i_+}, m]$ and $H_l$ is affine on $[x_{i_-}, x_{i_+}]$. Since $H_l$ is a piecewise affine, continuous function on $[a, m]$, we can define its right-hand continuous slope $h_l$. We thus get

$$h_l = h_{l-1} \qquad \text{on } \bar{J} = [a, x_{i_-}) \cup [x_{i_+}, m),$$
$$h_l = \bar{h}_{l-1}(J) \qquad \text{on } J = [x_{i_-}, x_{i_+}).$$

For all $l \geq 1$, the function $H_l$ is piecewise affine, continuous and $F \leq H_{l-1} \leq H_l \leq \widetilde{F}^m$ on $[a, m]$. The function $h_l$ is a cadlag step function on $[a, m)$ with discontinuity points in $X$. Moreover, using Lemma 3, we get $\|h_l - g\|([a, m)) \leq \|h_{l-1} - g\|([a, m))$. Therefore, $\|h_l - g\|([a, m)) \leq \|h_0 - g\|([a, m))$.

Within a finite number of iterations, we get $H_l = H_{l-1}$ and $h_l = h_{l-1}$. Let $l_0 + 1$ be the first step where it happens. Then $h_{l_0}$ is a nonincreasing function, so $H_{l_0}$ is concave. By definition of $\widetilde{F}^m$, we thus get $H_{l_0} \geq \widetilde{F}^m$ and thus, $H_{l_0} = \widetilde{F}^m$ on $[a, m]$. Therefore, $h_{l_0} = f^m$ on $[a, m)$. Finally,

$$\|f^m - g\|([a, m)) \leq \|h_0 - g\|([a, m)). \tag{12}$$

(2) Let us now see what is happening on $[m, b]$ in the case where $m \leq b$. The same holds if we replace $F$ by its left-hand continuous version $F^-$ on $(m, b]$, setting moreover $F^-(m) = F(m)$. Let us define $H_0$ on $[m, b]$ as the piecewise affine, continuous function such that on each $I_k$, $k > j$, $H_0$ is the greatest convex minorant of the restriction of $F^-$ to $I_k$. Let $h_0$ be its right-hand continuous slope. We iterate the symmetric rule: if $h_0$ is not nondecreasing, let $[x_{i-}, x_{i+}]$ be the first interval on which it happens. On this interval, we replace $h_0$ by its mean and $H_0$ by an affine function. We iterate this rule until we obtain a nondecreasing function $h_{l_0}$. We can check that $h_{l_0} = f^m$ on $[m, b)$ and that $\|h_{l_0} - g\|([m, b)) \leq \|h_0 - g\|([m, b))$. Thus,

$$\|f^m - g\|([m, b)) \leq \|h_0 - g\|([m, b)). \tag{13}$$

By summation of (12) and (13), we get

$$\|f^m - g\| \leq \|h_0 - g\|. \tag{14}$$

Now, by a straightforward decomposition, we get on each $I_k$

$$\|h_0 - g\|(I_k) \leq bh_0(I_k) + |\bar{h}_0(I_k) - \bar{g}(I_k)| l(I_k) + bg(I_k). \tag{15}$$

(3) Let $k \leq j$; by definition of $H_0$, we get $H_0(t_k) = F(t_k)$ and $H_0(t_{k-1}) = F(t_{k-1})$. Therefore,

$$|\bar{h}_0(I_k) - \bar{g}(I_k)| l(I_k) = |Z(t_k) - Z(t_{k-1})| \leq \sup_{t \in I_k} |Z(t) - Z(t_{k-1})|. \tag{16}$$

To compute the first term at the right-hand side of (15), we will use the following result (see Section 5.3 for the proof):

Lemma 4. *Let $h$ be a nonincreasing function on $J = [t_0, t_1]$. Let $H = \int_{t_0}^\cdot h(t) \, dt$. Then, using Notation 3,*

$$bh(J) = 2 \sup_{t \in J} (H(t) - H(t_0) - (t - t_0)\bar{h}(J)).$$

We apply Lemma 4 to the nonincreasing function $h_0$ on $I_k$. Since $H_0(t_{k-1}) = F(t_{k-1})$, $H_0(t_k) = F(t_k)$ and $H_0 \geq F$ on $I_k$, we get

$$bh_0(I_k) \geq 2 \sup_{t \in I_k} \left( F(t) - F(t_{k-1}) - \frac{t - t_{k-1}}{l(I_k)} (F(t_k) - F(t_{k-1})) \right). \tag{17}$$

But the restriction of $H_0$ to $I_k$ is a function whose slope can only change at the discontinuity points of $F$ where $H_0$ hits $F$, so that the supremum in Lemma 4 is achieved at such points. Let us call this set $Y$. We get

$$(18) \quad bh_0(I_k) = 2 \sup_{t \in Y} \left( F(t) - F(t_{k-1}) - \frac{t - t_{k-1}}{l(I_k)}(F(t_k) - F(t_{k-1})) \right).$$

Equations (17) and (18) yield

$$(19) \quad bh_0(I_k) = 2 \sup_{t \in I_k} \left( F(t) - F(t_{k-1}) - \frac{t - t_{k-1}}{l(I_k)}(F(t_k) - F(t_{k-1})) \right).$$

A last decomposition of the right-hand side of (19), using the concavity of $G$, gives

$$bh_0(I_k) \leq 4 \sup_{t \in I_k} |Z(t) - Z(t_{k-1})| + 2 \sup_{t \in I_k} (G(t) - G(t_{k-1}) - (t - t_{k-1})\bar{g}(I_k)).$$

Finally, applying Lemma 4 to $g$, we get

$$(20) \qquad bh_0(I_k) \leq 4 \sup_{t \in I_k} |Z(t) - Z(t_{k-1})| + bg(I_k).$$

Replacing (16) and (20) in (15) leads to

$$(21) \qquad \|h_0 - g\|(I_k) \leq 5 \sup_{t \in I_k} |Z(t) - Z(t_{k-1})| + 2bg(I_k).$$

(4) Let $k > j$; by definition of $H_0$, we get $H_0(t_{k-1}) = F^-(t_{k-1})$ and $H_0(t_k) = F^-(t_k)$. Setting $Z^- = F^- - G$, the second term on the right-hand side of (15) gives

$$(22) \quad \begin{aligned} &|\bar{h}_0(I_k) - \bar{g}(I_k)|l(I_k) \\ &= |Z^-(t_k) - Z^-(t_{k-1})| \\ &\leq |Z^-(t_k) - Z(t_{k-1}) + Z(t_{k-1}) - Z(t_{k-2}) + Z(t_{k-2}) - Z^-(t_{k-1})| \\ &\leq \sup_{t \in I_k} |Z(t) - Z(t_{k-1})| + 2 \sup_{t \in I_{k-1}} |Z(t) - Z(t_{k-2})|. \end{aligned}$$

To compute the first term on the right-hand side of (15), we apply Lemma 4 to the right-hand continuous nonincreasing function $-h_0$ on $I_k$,

$$bh_0(I_k) = 2 \sup_{t \in I_k} (H_0(t_{k-1}) - H_0(t) + (t - t_{k-1})\bar{h}_0(I_k)).$$

By construction, $H_0(t_{k-1}) = F^-(t_{k-1})$, $H_0(t_k) = F^-(t_k)$ and $H_0 \leq F^-$ on $I_k$. On the other hand, the slope of $H_0$ on $I_k$ can only change at the discontinuity points of $F^-$ such that $F^-(t) = H_0(t)$. Thus, using the same scheme as for (19), we get

$$bh_0(I_k) = 2 \sup_{t \in I_k} \left( F^-(t_{k-1}) - F^-(t) + \frac{t - t_{k-1}}{l(I_k)}(F^-(t_k) - F^-(t_{k-1})) \right).$$

Let $(y_j, y_{j+1})$ be the $j$th interval in $I_k$ on which $F^-$ is continuous. Then

$$bh_0(I_k) = 2 \sup_j \sup_{(y_j, y_{j+1})} \left( F^-(t_{k-1}) - F^-(t) + \frac{t - t_{k-1}}{l(I_k)}(F^-(t_k) - F^-(t_{k-1})) \right).$$

By continuity, the supremum on $(y_j, y_{j+1}]$ equals the supremum on the open interval $(y_j, y_{j+1})$, on which $F^- = F$. Therefore,

$$\begin{aligned} bh_0(I_k) &= 2 \sup_j \sup_{(y_j, y_{j+1})} \left( F^-(t_{k-1}) - F(t) + \frac{t - t_{k-1}}{l(I_k)}(F^-(t_k) - F^-(t_{k-1})) \right) \\ &= 2 \sup_{t \in I_k} \left( F^-(t_{k-1}) - F(t) + \frac{t - t_{k-1}}{l(I_k)}(F^-(t_k) - F^-(t_{k-1})) \right). \end{aligned}$$

A straightforward decomposition using Lemma 4 yields

$$bh_0(I_k) \le 2 \sup_{t \in I_k} |Z^-(t_{k-1}) - Z(t)| + 2|Z^-(t_k) - Z^-(t_{k-1})| + bg(I_k).$$

Finally,

$$(23) \quad bh_0(I_k) \le 4 \sup_{t \in I_k} |Z(t) - Z(t_{k-1})| + 8 \sup_{t \in I_{k-1}} |Z(t) - Z(t_{k-2})| + bg(I_k).$$

Replacing (22) and (23) in (15) leads to

$$\|h_0 - g\|(I_k) \le 5 \sup_{t \in I_k} |Z(t) - Z(t_{k-1})| + 10 \sup_{t \in I_{k-1}} |Z(t) - Z(t_{k-2})| + 2bg(I_k)$$

for all $k > j$. By (21) this inequality holds for all $k \in \mathcal{K}$, and by summation

$$(24) \quad \|h_0 - g\| \le \sum_{k \in \mathcal{K}} \left[ 15 \sup_{t \in I_k} |Z(t) - Z(t_{k-1})| + 2bg(I_k) \right].$$

Now let $p_\pi g$ denote the $\mathbb{L}_2$-orthogonal projection of $g$ on $\mathcal{H}_\pi$, that is,

$$(25) \quad p_\pi g(t) = \sum_{k \in \mathcal{K}} \bar{g}(I_k) \mathbb{1}_{[t_{k-1}, t_k)}(t) \qquad \text{for all } t \in I.$$

We get

$$\sum_{k \in \mathcal{K}} bg(I_k) = \|p_\pi g - g\|.$$

Setting $h \in \mathcal{H}_\pi$ and with $p_\pi h$ its $\mathbb{L}_2$-orthogonal projection on $\mathcal{H}_\pi$, we thus get $p_\pi h = h$ and then

$$(26) \quad \|p_\pi g - g\| \le 2\|g - h\|.$$

Therefore,

$$(27) \quad \sum_{k \in \mathcal{K}} bg(I_k) \le 2d(g, \mathcal{H}_\pi).$$

Substituting (27) in (24) completes the proof of Lemma 1, since we have (14), and $C' = 45$ works.

5.2. *Proof of Lemma* 2. The key arguments here are Marshall's lemma (see, e.g., [1]) and the following lemma (the proof of this lemma is omitted since it can be presented in the same way as that of Lemma 1 of [4]):

LEMMA 5. *Let* $F \in \mathcal{H}(I)$ *and let* $\widetilde{F}^r$ *and* $\widetilde{F}^s$ *be the U-shaped regularizations of* $F$ *on* $I$ *at* $r$ *and* $s$, *respectively, with* $r < s$. *Let* $f^r$ *and* $f^s$ *be their right-hand continuous slopes. Then*

$$\|f^r - f^s\| = 2\max\left\{\sup_{r \leq t \leq s}(F(t) - \widetilde{F}^r(t)), \sup_{r \leq t \leq s}(\widetilde{F}^s(t) - F(t))\right\}.$$

Let $m(F)$ be the point in $I$ such that $f = U_S^{m(F)}(F)$. By Lemma 5 we get

$$\|f^m - f\| \leq 2\max\left\{\sup_{t \in I}|F(t) - \widetilde{F}^m(t)|, \sup_{t \in I}|F(t) - \widetilde{F}^{m(F)}(t)|\right\}.$$

By definition of $m(F)$

$$\sup_{t \in I}|F(t) - \widetilde{F}^{m(F)}(t)| \leq \sup_{t \in I}|F(t) - \widetilde{F}^m(t)|,$$

and therefore

$$\|f^m - f\| \leq 2\sup_{t \in I}|F(t) - G(t)| + 2\sup_{t \in I}|G(t) - \widetilde{F}^m(t)|. \tag{28}$$

For the last term on the right-hand side of (28), we have

$$2\sup_{t \in I}|G(t) - \widetilde{F}^m(t)| = 2\max\left\{\sup_{t \leq m}|G(t) - \widetilde{F}^m(t)|, \sup_{t > m}|G(t) - \widetilde{F}^m(t)|\right\}, \tag{29}$$

so by Marshall's lemma,

$$2\sup_{t \in I}|G(t) - \widetilde{F}^m(t)| \leq 2\max\left\{\sup_{t \leq m}|G(t) - F(t)|, \sup_{t > m}|G(t) - F(t)|\right\} \tag{30}$$

$$\leq 2\sup_{t \in I}|F(t) - G(t)|.$$

Substituting (30) in (28) leads to Lemma 2.

5.3. *Proof of Lemma* 4. Let $u$ be the function defined on $J$ by $u = h - \bar{h}(J)$ and let $U$ be defined on $J$ by

$$U(t) = \int_{t_0}^t u(x)\,dx = H(t) - H(t_0) - (t - t_0)\bar{h}(J).$$

Then $u(t_0) \geq 0$ and $u(t_1) \leq 0$, so that there exists some $c \in J$ where $U$ achieves its maximum. Moreover, $u$ is nonnegative before $c$ and nonpositive after $c$. Since $U(t_1) = 0$, we thus get

$$bh(J) = \int_J |u(t)|\,dt = \int_{t_0}^c u(t)\,dt - \int_c^{t_1} u(t)\,dt = 2\sup_{t \in J}U(t),$$

which proves the lemma.

**6. Proof of Theorem 3.** In the sequel we adopt the convention and notation used in Section 5.

Let $\pi \in \Pi(I)$ and let $p_\pi g$ be the $\mathbb{L}_2$-orthogonal projection of $g$ on $\mathcal{H}_\pi$, defined by (25). To perform the right-hand side inequality of (5), we can write

$$\|\widehat{g}^\pi - g\| = \sum_{k \in \mathcal{K}} \int_{I_k} |g(t) - p_\pi g(t) + p_\pi g(t) - \widehat{g}^\pi(t)| \, dt$$

$$\leq \sum_{k \in \mathcal{K}} \int_{I_k} |g(t) - p_\pi g(t)| \, dt + \sum_{k \in \mathcal{K}} \int_{I_k} |p_\pi g(t) - \widehat{g}^\pi(t)| \, dt$$

$$\leq \|p_\pi g - g\| + \sum_{k \in \mathcal{K}} |Z(t_k) - Z(t_{k-1})|$$

$$\leq 2d(g, \mathcal{H}_\pi) + \sum_{k \in \mathcal{K}} \sup_{t \in I_k} |Z(t) - Z(t_{k-1})|.$$

The last control of $\|p_\pi g - g\|$ arises from (26). We get the result by a last obvious majorization ($C \geq 1$). Let us now prove the left-hand side inequality of (5). We get

$$\|\widehat{g}^\pi - g\| \geq \sum_{k \in \mathcal{K}} \left| \int_{I_k} (g(t) - p_\pi g(t) + p_\pi g(t) - \widehat{g}^\pi(t)) \, dt \right|$$

(31)
$$= \sum_{k \in \mathcal{K}} \left| \int_{I_k} (p_\pi g(t) - \widehat{g}^\pi(t)) \, dt \right|$$

$$= \sum_{k \in \mathcal{K}} |Z(t_k) - Z(t_{k-1})|.$$

On the other hand, by the triangle inequality,

(32)
$$\|\widehat{g}^\pi - g\| \geq \|g - p_\pi g\| - \|\widehat{g}^\pi - p_\pi g\|$$

$$= \|g - p_\pi g\| - \sum_{k \in \mathcal{K}} |Z(t_k) - Z(t_{k-1})|.$$

Multiplying (31) by $(CA + 4)$ and (33) by 4, and summing the so-obtained inequalities, we get, since $p_\pi g \in \mathcal{H}_\pi$,

$$(CA + 8)\|\widehat{g}_\pi - g\| \geq 4\|g - p_\pi g\| + CA \sum_{k \in \mathcal{K}} |Z(t_k) - Z(t_{k-1})|$$

$$\geq 4d(g, \mathcal{H}_\pi) + CA \sum_{k \in \mathcal{K}} |Z(t_k) - Z(t_{k-1})|.$$

Therefore, taking the expectations, we get that when condition (4) holds,

$$\mathcal{R}_Z(\pi) \leq 4d(g, \mathcal{H}_\pi) + CA \sum_{k \in \mathcal{K}} \mathbb{E}|Z(t_k) - Z(t_{k-1})|$$

$$\leq (CA + 8)\mathbb{E}\|\widehat{g}_\pi - g\|.$$

Relation (6) is straightforwardly derived from the last inequality and Theorem 2.

**7. Short proofs for Propositions 1–4.** Let us set $Z = \widehat{G} - G$. For each proposition, the proof follows the same scheme: to show (8) [resp. (10)], one needs to control the error term $\mathcal{R}_Z(\pi)$ in Theorem 2, for all $\pi \in \Pi(I)$. For this purpose we fix a partition $\pi$ in $\Pi_D(I)$, $D \le n$ (resp. $D \le T$), and we use the same notation as in the former section, setting moreover $G(I_k) = G(t_k) - G(t_{k-1})$. We then show that there exists some $C'$ such that

$$\sum_{k=1}^{D} \mathbb{E}\left( \sup_{t \in I_k} |Z(t) - Z(t_{k-1})| \right) \le C' \sqrt{\frac{D}{n}}, \tag{33}$$

respectively

$$\frac{1}{T} \sum_{k=1}^{D} \mathbb{E}\left( \sup_{t \in I_k} |Z(t) - Z(t_{k-1})| \right) \le 2\sqrt{\frac{D}{T}} \sqrt{\frac{G(T)}{T}}.$$

Therefore, we get the result applying Theorem 2 with $E = CC'$ (resp. $B = 2C$). Next, to show (9), one needs to check condition (4) in Theorem 3, for all $\pi \in \Pi(I)$.

PROOF OF PROPOSITION 1. To prove (33), let us call $F$ the common conditional distribution function of the $X_i$'s given that $X_i \in I_k$. Let $N = \sum_{i=1}^{n} \mathbb{1}_{X_i \in I_k}$ be the number of observations falling in $I_k$ and let $\widehat{F}_N$ be the empirical distribution of $N$ observations falling in $I_k$. We get for all $t \in I_k$,

$$F(t) = \frac{G(t) - G(t_{k-1})}{G(I_k)} \quad \text{and} \quad \widehat{F}_N(t) = \frac{n}{N}(\widehat{G}(t) - \widehat{G}(t_{k-1})).$$

Therefore,

$$\begin{aligned}
\mathbb{E}&\left( \sup_{t \in I_k} |\widehat{Z}(t) - \widehat{Z}(t_{k-1})| \right) \\
&\le \mathbb{E}\left( \sup_{t \in I_k} \frac{N}{n} |\widehat{F}_N(t) - F(t)| \right) + \mathbb{E}\left| \frac{N}{n} - G(I_k) \right|.
\end{aligned} \tag{34}$$

For the first term on the right-hand side of (34), an upper bound can be derived applying Massart's inequality [20] to $F$ on $I_k$,

$$\mathbb{P}\left( \sup_{t \in I_k} |\widehat{F}_N(t) - F(t)| > \lambda | N \right) \le 2e^{-2N\lambda^2} \qquad \forall \lambda > 0.$$

Integrating the latter inequality leads to

$$\mathbb{E}\left( \sup_{t \in I_k} \frac{N}{n} |\widehat{F}_N(t) - F(t)| \right) \le \sqrt{\pi/2} \, \mathbb{E}\left( \frac{\sqrt{N}}{n} \right). \tag{35}$$

A last control of (35) can be performed by the Cauchy–Schwarz's inequality, leading to

$$\mathbb{E}\left(\sup_{t \in I_k} \frac{N}{n} |\widehat{F}_N(t) - F(t)|\right) \leq \sqrt{\frac{\pi}{2}} \sqrt{\frac{G(I_k)}{n}}.$$

On the other hand, the second term on the right-hand side of (34) can be bounded by the Cauchy–Schwarz's inequality applied to $N \sim \mathcal{B}(n, G(I_k))$. We thus obtain

$$(36) \qquad \mathbb{E}\left(\sup_{t \in I_k} |Z(t) - Z(t_{k-1})|\right) \leq \left(1 + \sqrt{\frac{\pi}{2}}\right) \sqrt{\frac{G(I_k)}{n}}.$$

We then obtain (33), since

$$\sum_{k=1}^{D} \sqrt{\frac{G(I_k)}{n}} \leq \sqrt{\frac{D}{n}}.$$

To prove (9), one needs to sharpen the bound (36): actually, it gives the right order of magnitude of the supremum on each $I_k$ such that $G(I_k) \geq 1/n$, but it is too crude when $G(I_k) \leq 1/n$. Both $\widehat{G}$ and $G$ are monotone, so that for all $k$,

$$\mathbb{E}\left(\sup_{t \in I_k} |\widehat{Z}(t) - \widehat{Z}(t_{k-1})|\right) \leq 2G(I_k).$$

Combining this inequality with (36) yields

$$\mathbb{E}\left(\sup_{t \in I_k} |Z(t) - Z(t_{k-1})|\right) \leq \min\left\{2G(I_k), (1 + \sqrt{\pi/2})\sqrt{\frac{G(I_k)}{n}}\right\}.$$

On the other hand, by a lemma of Devroye and Györfi (see [12], page 25),

$$(37) \qquad \mathbb{E}|Z(t_k) - Z(t_{k-1})| \geq \min\left\{0.13 G(I_k), 0.36\sqrt{\frac{G(I_k)}{n}}\right\},$$

so there exists an absolute constant $A$ such that (4) holds. $\square$

PROOF OF PROPOSITION 2. Let $H$ and $B$ be the processes defined on $[0,1]$ by $H = \mathbb{E}(Z)$ and $B = Z - H$.

The $\varepsilon_i$'s are independent, so that applying the Cauchy–Schwarz inequality, we get

$$(38) \qquad \begin{aligned} \mathbb{E}\left(\sup_{t \in I_k} |B(t) - B(t_{k-1})|\right) &\leq \frac{1}{n}\mathbb{E}\left(\sum_{i=[nt_{k-1}]+1}^{[nt_k]} |\varepsilon_i|\right) \\ &\leq \frac{\sigma}{n}(\sqrt{n(t_k - t_{k-1})} + 1). \end{aligned}$$

On the other hand, $g$ is unimodal, so one can prove that

$$(39) \qquad \sup_{t \in I_k} |H(t) - H(t_{k-1})| \le \frac{6M}{n}.$$

Since

$$\sum_{i=1}^{D} \sqrt{\frac{l(I_k)}{n}} \le \sqrt{\frac{D}{n}},$$

relations (39) and (38) lead to (33), with $C' = 2(3M + \sigma)$.

To prove (9), let us only consider $k$ such that $[nt_k] \ne [nt_{k-1}]$ [otherwise (4) is trivial for all $A > 0$]. We use a Von Bahr–Esseen inequality (see, e.g., [27], page 858), which leads to

$$\mathbb{E}\Big(\sup_{t \in I_k} |B(t) - B(t_{k-1})|\Big) \le 8\mathbb{E}|B(t_k) - B(t_{k-1})|.$$

Using (39) and the preceding display, the triangle inequality gives

$$(40) \qquad 8\mathbb{E}(|Z(t_k) - Z(t_{k-1})|) \ge \mathbb{E}\Big(\sup_{t \in I_k} |B(t) - B(t_{k-1})|\Big) - 8\frac{6M}{n}.$$

Using Markov's inequality, we get

$$(41) \quad \mathbb{E}|Z(t_k) - Z(t_{k-1})| \ge \frac{\sigma\sqrt{2\pi}}{4n} \mathbb{P}\Big(|Z(t_k) - Z(t_{k-1})| \ge \frac{\sigma\sqrt{2\pi}}{4n}\Big) \ge \frac{\sqrt{2\pi}}{8n}.$$

The last inequality arises from the fact that $Z(t_k) - Z(t_{k-1})$ is a centered Gaussian variable whose variance is greater than $\sigma^2/n^2$.

Now, multiplying (41) by $432M/\sqrt{2\pi}$ and by summation with (40), there exists an $A'$ such that

$$\mathbb{E}\Big(\sup_{t \in I_k} |B(t) - B(t_{k-1})|\Big) \le A'\mathbb{E}|Z(t_k) - Z(t_{k-1})|.$$

The triangle inequality and relations (39) and (41) lead to the fact that there exists an $A$ such that condition (4) holds. $\square$

PROOF OF PROPOSITION 3. Let us set for all $t \in I$ $g^*(t) = \mathbb{1}_{X_{(n)} \ge t} g(t)$ and $G^*(t) = \int_0^t g^*(s)\,ds$, where $X_{(n)}$ is the $n$th order statistic of the sample. Setting $Z^* = \widehat{G} - G^*$, we get

$$(42) \begin{aligned} &\sum_{k=1}^{D} \mathbb{E}\Big(\sup_{t \in I_k} |Z(t) - (t_{k-1})|\Big) \\ &\le \sum_{k=1}^{D} \mathbb{E}\Big(\sup_{t \in I_k} |Z^*(t) - Z^*(t_{k-1})|\Big) + \mathbb{E}\Big(\int_0^c \mathbb{1}_{X_{(n)} < s} g(s)\,ds\Big). \end{aligned}$$

The process $(Z^*(t))_{t \geq 0}$ is a square integrable mean zero martingale. Its predictable variation process is given by (see, e.g., [27], Theorem 2, page 312)

$$\langle Z^* \rangle = \int_0^\cdot \frac{\mathbb{1}_{s \leq X_{(n)}}}{n(1 - \widehat{L}^-(t))} g(s) \, ds, \tag{43}$$

where $\widehat{L}^-$ is the left-hand continuous version of the empirical distribution function of the $X_i$'s. Using Doob's inequality, relation (43) combined with the Cauchy–Schwarz inequality yields for all $k$,

$$\mathbb{E}\left(\sup_{t \in I_k} |Z^*(t) - Z^*(t_{k-1})|\right) \leq 2\sqrt{\mathbb{E}(\langle Z^* \rangle(t_k) - \langle Z^* \rangle(t_{k-1}))}$$

$$\leq \frac{2}{\sqrt{n}} \sqrt{\int_{I_k} \mathbb{E}\left(\sup_{s \leq X_{(n)}} \frac{1 - L(s)}{1 - \widehat{L}^-(s)}\right) \frac{g(s)}{1 - L(s)} \, ds}.$$

Setting $H$ as the common distribution function of the $U_i$'s [we get $1 - L = (1 - F)(1 - H)$], we thus apply Gill's inequality (see [27]), which leads to

$$\mathbb{E}\left(\sup_{t \in I_k} |Z^*(t) - Z^*(t_{k-1})|\right) \leq \sqrt{\frac{12}{n}} \frac{1}{\sqrt{1 - H(c)}} \sqrt{\frac{1}{1 - F(t_k)} - \frac{1}{1 - F(t_{k-1})}}.$$

Finally,

$$\sum_{k=1}^{D} \mathbb{E}\left(\sup_{t \in I_k} |Z^*(t) - Z^*(t_{k-1})|\right) \leq \sqrt{\frac{D}{n}} \frac{\sqrt{12}}{\sqrt{1 - L(c)}}. \tag{44}$$

For the last term on the right-hand side of (42), we use the relations

$$\mathbb{P}(X_{(n)} \leq s) = (1 - (1 - H(s))(1 - F(s)))^n \quad \text{and} \quad \left(1 - \frac{s}{n}\right)^n \leq e^{-s}$$

$$\text{for all } s \leq n.$$

Simple calculations yield

$$\mathbb{E}\left(\int_0^c \mathbb{1}_{X_{(n)} < s} g(s) \, ds\right) \leq \sqrt{\frac{\pi}{n}} \frac{1}{\sqrt{1 - L(c)}}. \tag{45}$$

Combining (44) and (45) in (42) yields (33). □

PROOF OF PROPOSITION 4. The process $(Z(t))_{t \geq 0}$ is a square integrable mean zero martingale. Since $N(t_k) - N(t_{k-1}) \sim \mathcal{P}(G(I_k))$, applying the Doob and Cauchy–Schwarz inequalities thus leads to

$$\mathbb{E}\left(\sup_{t \in I_k} |Z(t) - Z(t_{k-1})|\right) \leq 2\sqrt{G(I_k)},$$

so that

$$\frac{1}{T}\sum_{k=1}^{D}\mathbb{E}\left(\sup_{t\in I_k}|Z(t)-Z(t_{k-1})|\right)\leq 2\sqrt{\frac{D}{T}}\sqrt{\frac{G(T)}{T}}.$$

This proves (33).

To prove (9), let us set for all $x\geq 0$ $M(x)=N_1(x)-x$, where $(N_1(x))_{x\geq 0}$ is the Poisson process with mean function $x$. Recall that $X_j=G(T_j)$ is the $j$th occurrence time of $(N_1(x))_{x\geq 0}$. Then let us set $J=[x_{k-1},x_k]$, where $x_k=G(t_k)$ and $x_{k-1}=G(t_{k-1})$ and let $(a_i)_{0\leq i\leq m}$ be the sequence of endpoints of a uniform partition of $J$. Since $N_1$ has independent increments, the $m$ random variables $M(a_i)-M(a_{i-1})$ are integrable i.i.d. mean zero variables and we can apply a Von Bahr–Esseen inequality (see [27], page 858),

$$\mathbb{E}\left(\max_{1\leq i\leq m}|M(a_i)-M(x_{k-1})|\right)\leq 8\mathbb{E}|M(x_k)-M(x_{k-1})|.$$

Moreover, $N_1$ is a cadlag process. Therefore, this inequality holds on the whole interval $J$, when the partition's step tends toward zero. We thus get

$$\mathbb{E}\left(\sup_{x\in J}|M(x)-M(x_{k-1})|\right)\leq 8\mathbb{E}(|M(x_k)-M(x_{k-1})|)$$

and then

$$\mathbb{E}\left(\sup_{t\in I_k}|Z(t)-Z(t_{k-1})|\right)\leq 8\mathbb{E}(|Z(t_k)-Z(t_{k-1})|). \qquad \square$$

**Acknowledgments.** This paper is an extension of my Ph.D. thesis supervised by P. Massart, Professor at the University Paris-Sud. I thank him for his precious help. I also thank C. Durot for her useful comments which have helped to considerably improve the quality of this paper.

## REFERENCES

[1] Barlow, R. E., Bartholomew, D. J., Bremner, J. M. and Brunk, H. D. (1972). *Statistical Inference under Order Restrictions*. Wiley, New York. MR326887

[2] Barlow, R. E., Proschan, F. and Scheuer, E. M. (1972). A system debugging model. In *Reliability Growth Symposium, Interim Note 22 Aberdeen Proving Ground* 46–65. U.S. Army Materiel Systems Analysis Agency, Washington, DC. MR368368

[3] Bartoszyński, R., Brown, B. W., McBride, C. M. and Thompson, J. R. (1981). Some nonparametric techniques for estimating the intensity function of a cancer related nonstationary Poisson process. *Ann. Statist.* **9** 1050–1060. MR628760

[4] Birgé, L. (1987). Robust estimation of unimodal densities. Technical report, Univ. Paris X–Nanterre.

[5] Birgé, L. (1987). On the risk of histograms for estimating decreasing densities. *Ann. Statist.* **15** 1013–1022. MR902242

[6] Birgé, L. (1989). The Grenander estimator: A nonasymptotic approach. *Ann. Statist.* **17** 1532–1549. MR1026298

[7] BIRGÉ, L. (1997). Estimation of unimodal densities without smoothness assumptions. *Ann. Statist.* **25** 970–981. MR1447736

[8] BRUNK, H. D. (1970). Estimation of isotonic regression. In *Nonparametric Techniques in Statistical Inference* (M. L. Puri, ed.) 177–197. Cambridge Univ. Press, London. MR277070

[9] CENCOV, N. N. (1962). Evaluation of an unknown distribution density from observations. *Soviet Math. Dokl.* **3** 1559–1562.

[10] CLEVENSON, M. L. and ZIDEK, J. V. (1977). Bayes linear estimators of the intensity function of the nonstationary Poisson process. *J. Amer. Statist. Assoc.* **72** 112–120. MR461826

[11] CURIONI, M. (1977). Estimation de la densité des processus de Poisson non homogènes. Ph.D. dissertation, Univ. Pierre et Marie Curie–Paris VI.

[12] DEVROYE, L. and GYÖRFI, L. (1985). *Nonparametric Density Estimation: The $L_1$ View.* Wiley, New York. MR780746

[13] DEVORE, R. A. and LORENTZ, G. G. (1993). *Constructive Approximation.* Springer, Berlin. MR1261635

[14] DUROT, C. (2002). Sharp asymptotics for isotonic regression. *Probab. Theory Related Fields* **122** 222–240. MR1894068

[15] GRENANDER, U. (1956). On the theory of mortality measurement. II. *Skand.-Aktuarietidskr.* **39** 125–153. MR93415

[16] GROENEBOOM, P. (1985). Estimating a monotone density. In *Proc. Berkeley Conference in Honor of Jerzy Neyman and Jack Kiefer* (L. M. Le Cam and R. A. Olshen, eds.) **2** 539–555. Wadsworth, Monterey, CA. MR822052

[17] GROENEBOOM, P., HOOGHIEMSTRA, G. and LOPUHAÄ, H. P. (1999). Asymptotic normality of the $L_1$ error of the Grenander estimator. *Ann. Statist.* **27** 1316–1347. MR1740109

[18] KAPLAN, E. L. and MEIER, P. (1958). Nonparametric estimation from incomplete observations. *J. Amer. Statist. Assoc.* **53** 457–481. MR93867

[19] MARRON, J. S. and PADGETT, W. J. (1987). Asymptotically optimal bandwidth selection for kernel density estimators from randomly right-censored samples. *Ann. Statist.* **15** 1520–1535. MR913571

[20] MASSART, P. (1990). The tight constant in the Dvoretzky–Kiefer–Wolfowitz inequality. *Ann. Probab.* **18** 1269–1283. MR1062069

[21] MÜLLER, H.-G. and WANG, J.-L. (1990). Locally adaptive hazard smoothing. *Probab. Theory Related Fields* **85** 523–538. MR1061943

[22] MÜLLER, H.-G. and WANG, J.-L. (1994). Hazard rate estimation under random censoring with varying kernels and bandwidths. *Biometrics* **50** 61–76. MR1279435

[23] NELSON, W. (1972). Theory and applications of hazard plotting for censured failure data. *Technometrics* **14** 945–966.

[24] PARZEN, E. (1962). On estimation of a probability density function and mode. *Ann. Math. Statist.* **33** 1065–1076. MR143282

[25] PRIESTLEY, M. B. and CHAO, M. T. (1972). Non-parametric function fitting. *J. Roy. Statist. Soc. Ser. B* **34** 385–392. MR331616

[26] RUDEMO, M. (1982). Empirical choice of histograms and kernel density estimators. *Scand. J. Statist.* **9** 65–78. MR668683

[27] SHORACK, G. R. and WELLNER, J. A. (1986). *Empirical Processes with Applications to Statistics.* Wiley, New York. MR838963

[28] SINGPURWALLA, N. D. and WONG, M. Y. (1983). Estimation of the failure rate. A survey of nonparametric methods. I. Non-Bayesian methods. *Comm. Statist. Theory Methods* **12** 559–588. MR696809

[29] Stone, C. J. (1984). An asymptotically optimal window selection rule for kernel density estimates. *Ann. Statist.* **12** 1285–1297. MR760688

[30] Tanner, M. A. and Wong, W. H. (1983). The estimation of the hazard function from randomly censored data by the kernel method. *Ann. Statist.* **11** 989–993. MR707949

[31] Wang, Y. (1995). The $L_1$ theory of estimation of monotone and unimodal densities. *J. Nonparametr. Statist.* **4** 249–261. MR1366772

[32] Wegman, E. J. (1970). Maximum likelihood estimation of a unimodal density function. *Ann. Math. Statist.* **41** 457–471. MR254995

[33] Yandell, B. S. (1983). Nonparametric inference for rates with censored survival data. *Ann. Statist.* **11** 1119–1135. MR720258

14 rue de Bretagne
75003 Paris
France
e-mail: laurence.reboul@wanadoo.fr